\documentclass[12pt,a4paper]{article}
\usepackage{amsmath,amssymb,amsfonts}
\def\Bbb{\mathbb}

\title{\bf On the values of Dedekind sums}

\author{Kurt Girstmair}
\date{}

\makeatletter
\let\@@maketitle=\maketitle
\def\maketitle{\def\thispagestyle##1{\relax}\@@maketitle}
\makeatother
\newtheorem{theorem}{Theorem}
\newtheorem{prop}{Proposition}

\def\BE{\begin{equation}}
\def\EE{\end{equation}}
\def\BD{\begin{displaymath}}
\def\ED{\end{displaymath}}
\def\BA{\begin{array}}
\def\EA{\end{array}}
\def\BEA{\begin{eqnarray*}}
\def\EEA{\end{eqnarray*}}
\def\BI{\bibitem}

\def\N{\Bbb N}
\def\Z{\Bbb Z}

\def\R{\Bbb R}

\def\phi{\varphi}

\def\MB{\mbox}
\def\LD{\ldots}
\def\OV{\overline}

\def\DIV{\,|\,}
\def\NDIV{\, \nmid \,}

\def\BQ{``}
\def\EQ{'' }

\def\MN{\medskip\noindent}
\def\STOP{\hfill$\Box$}

\def\JS#1#2{ \left( \frac{#1}{#2} \right) }
\def\DED{Dedekind }

\begin{document}
\maketitle

\begin{abstract}

\noindent
Let $s(a,b)$ denote the classical Dedekind sum and $S(a,b)=12s(a,b)$. For a given
denominator $q\in \N$, we study the numerators $k\in\Z$ of the values $k/q$, $(k,q)=1$, of Dedekind sums $S(a,b)$.
Our main result says that if $k$ is such a numerator, then the whole residue class of $k$ modulo $(q^2-1)q$ consists of numerators of this kind. This fact reduces the task of finding
all possible numerators $k$ to that of finding representatives for finitely many residue classes modulo $(q^2-1)q$.
By means of the proof of this result we have determined all possible numerators $k$ for $2\le q\le 60$, the case $q=1$ being trivial. The result of this search
suggests a conjecture about all possible values $k/q$, $(k,q)=1$, of Dedekind sums $S(a,b)$ for an arbitrary $q\in\N$.

\end{abstract}

\section*{1. Introduction and main results}

Let $a$ be an integer, $b$ a natural number, and $(a,b)=1$. The classical \DED sum $s(a,b)$ is defined by
\BD
   s(a,b)=\sum_{k=1}^{b} ((k/b))((ak/b)).
\ED
Here
\BD
  ((x))=\begin{cases}
                 x-\lfloor x\rfloor-1/2 & \MB{ if } x\in\R\smallsetminus \Z; \\
                 0 & \MB{ if } x\in \Z
        \end{cases}
\ED
(see \cite[p. 1]{RaGr}).
It is often more convenient to work with
\BD
 S(a,b)=12s(a,b)
\ED instead. We call $S(a,b)$ a {\em normalized} \DED sum.

The values of normalized \DED sums are rational numbers, whose properties have been studied by many authors.
For instance, it is known that they are dense in the field $\R$ (see \cite{Hi}); and zones of large
and small values of normalized \DED sums have been described by the present author (see \cite{Gi2}). It is unknown,
however, which rational numbers occur as the values of normalized \DED sums. It seems that there is not even a
conjecture in this direction. In this paper we present a conjecture of this kind that is based on empirical
data.

Our first main result says that the numerators of  normalized \DED sums with a given denominator $q$ do not appear isolated, but form
certain residue classes.

\begin{theorem} 
\label{t1}
Let $q$ be a natural number and $k$ an integer with $(k,q)=1$. If $k/q$ is the value of a normalized \DED sum, then all numbers
\BD
 k'/q \:\MB{ with }\: k'\in \Z,\enspace k'\equiv k\mod (q^2-1)q,
\ED
are values of normalized \DED sums.
\end{theorem} 

\MN
The case $q=1$ of Theorem \ref{t1} is trivial, since
$(q^2-1)q=0$ and 0 is the only integer value of a normalized \DED sum (see the first remark following the proof of Theorem \ref{t3} below).

Theorem 1 reduces the task of finding all values $k/q$, $(k,q)=1$, of normalized \DED sums to finitely many cases for each $q$.
One just has to find a representative $k$ of each residue class modulo $(q^2-1)q$ in question.

The {\em proof} of this theorem provides the tools for finding these representatives in particular cases.
In this way we have determined, for $2\le q \le 60$, all possible numerators $k$ of normalized \DED sums $S(a,b)=k/q$, $(k,q)=1$. The result is as follows.

\begin{theorem} 
\label{t2}
Let $q$ be a natural number, $2\le q\le 60$, and $k$ an integer with $(k,q)=1$. Then $k/q$ is the value of a normalized \DED sum, if, and only if,
$k$ satisfies the following conditions.
\\ {\rm(a)} If $3\NDIV q$, then $k\equiv 0\mod 3$.
\\ {\rm(b)} If $2\NDIV q$, then
\BD
   k\equiv\begin{cases}
           2\mod 4, \enspace \MB{ if } q\equiv 3\mod 4;\\
           0\mod 8, \enspace \MB{ if } q \MB{ is a square;}\\
           0\mod 4, \enspace \MB{ otherwise.}
           \end{cases}
\ED
\end{theorem} 

\MN
We will show that the conditions (a) and (b) necessarily hold if an arbitrary rational number $k/q$, $(k,q)=1$, is the value of a normalized \DED sum (see Propositions \ref{p1} and \ref{p2} below).
By Theorem \ref{t2}, these conditions are sufficient for $2\le q\le 60$. This fact leads to the following conjecture.

\MN
{\em Conjecture.} Let $q\ge 2$ be a natural number, $k\in\Z$, $(k,q)=1$. Then $k/q$ is the value of a normalized \DED sum if, and only if,
$k$ satisfies the conditions (a) and (b) of Theorem \ref{t2}.

\section*{2. The details}

The basis of the present work is the following theorem.

\begin{theorem} 
\label{t3}
Let $b, q$ be natural numbers, $a$ an integer, $(a,b)=1$. Then $S(a,b)$ takes the form $k/q$ for some $k\in\Z$, $(k,q)=1$ if,
and only if, $b$ has the form
\BE
\label{2.2}
 b=\frac{q(a^2+1)}t,
\EE
where $t$ is a natural number, $(t,q)=1$.
$($Note that the condition $(t,q)=1$ implies $t\DIV a^2+1$.$)$

\end{theorem} 

\MN
{\em Proof.} It suffices to assume $a\in\N$ since $S(-a,b)=-S(a,b)$ and $a^2+1=(-a)^2+1$.
The reciprocity law for normalized \DED sums (see \cite[p. 5]{RaGr}) gives
\BD
 abS(a,b)=-abS(b,a)+ a^2+b^2+1-3ab.
\ED
By \cite[p. 27, Th. 2]{RaGr}, $aS(b,a)$ is an integer, so $abS(a,b)\equiv a^2+1 \mod b$ or
\BD
  aS(a,b)\equiv \frac{a^2+1}b \mod \Z.
\ED
This means that $q=b/(b,a^2+1)$ is the smallest possible denominator of $aS(a,b)$. Since $bS(a,b)\in\Z$ and
$(a,b)=1$, it is also the smallest possible denominator of $S(a,b)$. Hence we may write
$b=q(a^2+1)/t$ with some natural number $t$ dividing $a^2+1$. It is easy to see that $(q,t)=1$, since, otherwise,
we would obtain a smaller denominator than $q$. Conversely, if $b$ has the form (\ref{2.2}) with $(t,q)=1$,
one sees that $(a^2+1)/t=(b, a^2+1)$. Indeed, if $(a^2+1)/t\cdot r$ divides both $b$ and $a^2+1$, then $r$ divides both $q$ and $t$, and so
$r=1$. Accordingly, $S(a,b)$ has the smallest possible denominator $q=bt/(a^2+1)$.
\STOP

\MN
{\em Remarks.} 1. For $q=1$, Theorem \ref{t3} shows that $b$ takes the form $b=(a^2+1)/t$. Hence $a^2\equiv -1\mod b$.
This, however, implies that $S(a,b)=0$, since $-a$ is a multiplicative inverse of $a$ modulo $b$, and so $S(a,b)=S(-a,b)=-S(a,b)$.
In other words, $0$ is the only integer that is the value of a normalized \DED sum. This fact
is well-known (see \cite[p. 36]{RaGr}).

2. A weaker form of Theorem \ref{t3} can be found in the author's paper \cite{Gi3}.

\MN
In what follows we assume that $S(a,b)=k/q$, $k\in\Z, q\in\N$, $(k,q)=1$. By Theorem \ref{t3} this is the same as saying
$b=q(a^2+1)/t$ for some $t\in\N$ with $(t,q)=1$. The following theorem expresses $k/q$ in terms of two normalized \DED sums.

\begin{theorem} 
\label{t4}
Let $q,t$ be natural numbers $(t,q)=1$, $a$ an integer, $(a,q)=1$, and $b=q(a^2+1)/t$.
Then
\BE
\label{2.4}
 S(a,b)=\frac{(q^2-1)a}{tq}-S(aq,t)+S(at^*,q),
\EE
where $t^*$ is an integer satisfying $tt^*\equiv 1 \mod q$.
\end{theorem} 

\MN
{\em Proof.} First we assume $a\in \N$. Then the reciprocity law gives
\BE
\label{2.6}
   S(a,b)=-S(b,a)+\frac ab+\frac ba+\frac 1{ab}-3.
\EE
We apply the three-term relation for normalized \DED sums to $S(b,a)$ (see \cite[Lemma 1]{Gi}).
It gives
\BE
\label{2.8}
S(b,a)=S(aq,t)+S(r,q)+\frac a{tq}+\frac t{aq}+\frac q{at}-3,
\EE
where $r$ is defined as follows: Let $j,l$ be integers such that
$-aqj+tl=1$. Then $r=-al+bj$. However, we only need the residue class of $r$ mod $q$.
Obviously, we have $l\equiv t^*\mod q$ and $r\equiv -al\mod q$. So we may choose $r=-at^*$.
Now equations (\ref{2.6}) and (\ref{2.8}), in combination with $b=q(a^2+1)/t$,
yield (\ref{2.4}).

In the case $a<0$ we write $S(a,b)=-S(|a|,b)$ and expand $S(|a|,b)$ as above. In the end,
we obtain (\ref{2.4}) again.
\STOP

\MN
{\em Proof of Theorem \ref{t1}}. Let $k\in\Z$, $q\in \N$, $(k,q)=1$, and suppose $k/q=S(a,b)$.
By Theorem \ref{t3}, $b=q(a^2+1)/t$ with $t\in\N$, $(t,q)=1$. By Theorem \ref{t4},
$S(a,b)$ is given by (\ref{2.4}). Let $m\in \Z$ and $a'=a+mtq$, $b'=q(a'^2+1)/t$ (observe that $a'^2+1\equiv a^2+1\equiv 0\mod t$).
Then (\ref{2.4}) gives
\BD
  S(a',b')=\frac{(q^2-1)a'}{tq}-S(a'q,t)+S(a't^*,q).
\ED
Here, however, $a'/(tq)= a/(tq)+m$; moreover, $S(a'q,t)=S(aq,t)$, since $a'\equiv a\mod t$, and $S(a't^*,q)=S(at^*,q)$,
since $a'\equiv a \mod q$. Altogether, we have
\BD
 S(a',b')=(q^2-1)m+S(a,b) =(q^2-1)m+\frac kq.
\ED
If we write $(q^2-1)m+k/q$ in the form $((q^2-1)qm+k)/q$, we see that $k'/q$ is the value of a normalized \DED sum for each number $k'\equiv k\mod (q^2-1)q$.
\STOP

\MN
The next two propositions show that the conditions (a) and (b) of Theorem \ref{t2} are necessary for $k/q$, $(k,q)=1$, being the value of a normalized \DED sum.

\begin{prop} 
\label{p1}
Let $q\in \N$ and $k\in\Z$, $(k,q)=1$ be such that $k/q$ is the value of a normalized \DED sum.
If $3\NDIV q$, then $k\equiv 0\mod 3$.
If $q\equiv 3\mod 4$, then $k\equiv 2\mod 4$. If $q\equiv 1\mod 4$, then
$k\equiv 0\mod 4$.

\end{prop} 

\MN
{\em Proof.} Suppose that $k/q=S(a,b)$ with $b=q(a^2+1)/t$ and $(t,q)=1$.
If $3\NDIV q$, we have $3\NDIV b$, since $a^2+1/t\not\equiv 0\mod 3$ (note that the prime divisors $p$ of $a^2+1$ are either $p=2$ or satisfy $p\equiv 1\mod 4$, $-1$ being a quadratic residue mod $p$).
In this case
$bS(a,b)\equiv 0 \mod 3$ (see \cite[p. 27, Th. 2]{RaGr}) and so $k\equiv 0$ mod 3.

Now suppose that $q$ is odd. We first consider the case that $b$ is odd. A classical result of \DED says
\BE
\label{2.10}
   bS(a,b)\equiv b+1-2\JS ab \mod 8,
\EE
where $\JS ab$ is the Jacobi symbol (see \cite[p. 34, formula (42)]{RaGr}). Since $(a^2+1)/t$ is odd, it must be $\equiv 1\mod 4$, by the above argument. Accordingly, $b\equiv q\mod 4$.
If $q$ (and, thus, $b$) is $\equiv 3$ mod $4$, then (\ref{2.10}) gives $bS(a,b)\equiv 2\mod 4$ and so $k\equiv 2\mod 4$. In the case $q\equiv 1\mod 4$ we obtain, in the same way,
$k\equiv 0 \mod 4$.

Let $q$ be odd and $b$ even. Then $(a^2+1)/t$ is even. Further, $t$ must be odd, since otherwise the common divisor $2$ of $a^2+1$ and $t$ disappears from $b$. Put $a'=a+tq$ and $b'=q(a'^2+1)/t$.
From $a'^2+1=a^2+1+2atq+t^2q^2$ we see that $(a'^2+1)/t$ is odd, and so $b'$ is odd. Theorem \ref{t4} shows that $S(a',b')=k'/q$ with $k'=k+(q^2-1)q$. However, $q^2-1\equiv 0\mod 8$,
and so $k'\equiv k\mod 8$. From the above we have $k'\equiv 2\mod 4$ if $q\equiv 3\mod 4$, and $k'\equiv 0\mod 4$ if $q\equiv 1\mod 4$. The respective congruences also hold for
$k$ instead of $k'$.
\STOP

\begin{prop} 
\label{p2}
Let $q$ be an odd square and $k\in\Z$, $(k,q)=1$, be such that $k/q$ is the value of a normalized \DED sum. Then $k\equiv 0\mod 8$.

\end{prop} 

\MN
{\em Proof.} Let $k/q=S(a,b)$ with $b=q(a^2+1)/t$, $(t,q)=1$.
\\ {\em Case 1}: Suppose that $b$ is odd and $a$ is even. Then $t$ is odd. We write $a=2^ja'$ with an odd number $a'$. In view of (\ref{2.10}),
we have to evaluate $\JS ab$. Obviously,
\BE
\label{2.11}
  \JS ab={\JS 2b}^j\JS {a'}b.
\EE
Since $b\equiv 1\mod 4$,
\BD
 \JS {a'}b=\JS b{a'}=\JS{q(a^2+1)/t}{a'}=\JS q{a'}\JS t{a'},
\ED
because
\BD
  \JS{(a^2+1)/t}{a'}\JS t{a'}=\JS{a^2+1}{a'}=1.
\ED
However, $\JS q{a'}=1$, and so
\BE
\label{2.12}
  \JS{a'}b=\JS t{a'}=\JS{a'}t.
\EE
On the other hand,
\BD
 \JS 2b=\JS 2q\JS 2{a^2+1}\JS 2t.
\ED
Here $\JS 2q=1$. It is easy to see that
\BE
\label{2.14}
  \JS 2{a^2+1}= \begin{cases}
                 -1,        & \MB{ if } j=1;\\
                 \enspace 1,& \MB{ if } j\ge 2.
                \end{cases}
\EE
From (\ref{2.11}), (\ref{2.12}) and (\ref{2.14}) we conclude
\BE
\label{2.16}
  \JS ab= \begin{cases}
                 -\JS at,& \MB{ if } j=1;\\
                 \rule{0mm}{5mm}
                 \enspace\JS at ,& \MB{ if } j\ge 2.
                \end{cases}
\EE
In order to evaluate $\JS at$, we consider an arbitrary prime divisor $p$ of $t$. Since $a^2\equiv -1\mod p$,
$\OV a$ is an element of order 4 in the prime residue group $(\Z/\Z p)^{\times}$. Accordingly, $\OV a$ is a
square if, and only if, 8 divides the order of this group, i.e., $p\equiv 1\mod 8$. This gives
\BE
\label{2.18}
\JS ap= \begin{cases}
                 \enspace 1,        & \MB{ if } p\equiv 1\mod 8;\\
                  -1,& \MB{ if } p\equiv 5\mod 8.
                \end{cases}
\EE
If we write $t=p_1\cdots p_r$ with (not necessarily distinct) prime factors $p_i$ and apply (\ref{2.18}) to each of these factors, we see
that $\JS at=1$ if, and only if, $t\equiv 1\mod 8$. In view of (\ref{2.16}), we obtain that $\JS ab=1$ if, and only if, $j=1$ and $t\equiv 5\mod 8$
or $j\ge 2$ and $t\equiv 1\mod 8$. This is equivalent to $(a^2+1)/t\equiv 1\mod 8$. Since $q\equiv 1\mod 8$, we have
\BD
 \JS ab= \begin{cases}
                 \enspace 1,        & \MB{ if } b\equiv 1\mod 8;\\
                  -1,               & \MB{ if } b\equiv 5\mod 8.
                \end{cases}
\ED
Accordingly, (\ref{2.10}) gives $bS(a,b)\equiv 0\mod 8$ and $k\equiv 0\mod 8$.

{\em Case 2}: $b$ and $a$ are both odd. Then $t$ is even, more precisely, $t=2t'$ with $t'\equiv 1\mod 4$.
We have
\BE
\label{2.20}
  \JS ab=\JS ba=\JS qa\JS{a^2+1}a\JS ta=\JS ta=\JS 2a\JS{t'}a=\JS 2a\JS a{t'}.
\EE
Now
\BD
 \JS 2a=\begin{cases}
                 \enspace 1,        & \MB{ if } a^2\equiv 1\mod 16;\\
                  -1,& \MB{ if } a^2\equiv 9\mod 16.
                \end{cases}
\ED
Moreover, $(a^2+1)/2\equiv 1\mod 8$ if, and only if, $a^2\equiv 1\mod 16$, i.e., $\JS 2a=1$.
As in case 1, one shows that $\JS a{t'}=1$ if, and only if, $t'\equiv 1\mod 8$.
On observing
\BD
b=q\frac{a^2+1}2/t'
\ED
 and $q\equiv 1\mod 8$, we have
$b\equiv 1\mod 8$ if, and only if, $(a^2+1)/2\equiv 1\mod 8$ and $t'\equiv 1\mod 8$ or $(a^2+1)/2\equiv 5\mod 8$ and $t'\equiv 5\mod 8$. By
the above, however, this is equivalent to $\JS 2a=1$ and $\JS a{t'}=1$ or $\JS 2a=-1$ and $\JS a{t'}=-1$.
In view of (\ref{2.20}), we obtain that $b\equiv 1\mod 8$ if, and only if, $\JS ab=1$. Again, (\ref{2.10}) shows $k\equiv 0\mod 8$.

{\em Case 3}: $b$ is even. Then $t$ must be odd. We reduce this case to the case \BQ $b$ odd\EQ just as we did in the proof of Proposition \ref{p1}:
Put $a'=a+tq$, $b'=q(a'^2+1)/t$. Then $b'$
is odd and $S(a',b')=k'/q$ with $k'\equiv k\mod 8$. By the cases 1 and 2, $k'\equiv 0\mod 8$, hence $k\equiv 0\mod 8$.
\STOP

\MN
The {\em proof of Theorem \ref{t2}} is delivered by the computer, of course.
We describe the search procedure that gave us this theorem. In view of Theorem \ref{t4}, we worked with a list of pairs $(t,a_1)$, $t\ge 1$, $a_1\ge 0$, such that $(t,q)=1$, $0\le a_1\le t/2$ and $a_1^2+1\equiv 0\mod t$.
Note that $t$ may occur in this list more than once
if it is the product of two ore more distinct primes $\equiv 1\mod 4$. To any such pair we determined the numbers $a=a_1+tj$, $j=0,\LD,q-1$, such that $(a,q)=1$. In this way $a$ yields all residues
$a_2$, $0\le a_2\le q-1$, $(a_2,q)=1$ modulo $q$ (and the corresponding \DED sums $S(a_2t^*,q)$ in the sense of Theorem \ref{t4}). Since $a\equiv a_1\mod t$, $t$ divides $a^2+1$. Therefore, $b=q(a^2+1)/t\in\N$.
Of course, $(a,b)=1$.

We computed $S(a,b)$ and determined the residue class $\OV k$ modulo $(q^2-1)q$ of the numerator of $S(a,b)$ by a representative $k$. In the case
$q=2$ there is only the residue class $\OV 3$ modulo $(q^2-1)q$ satisfying the conditions (a) and (b) of Theorem \ref{t2}.
If $q\ge 3$, these residue classes take the form $\OV{\pm k}$, $1\le k< (q^2-1)q/2$. Observe that it suffices to find one of the numbers $\pm k$, since the other one arises from $S(-a,b)=-S(a,b)$.

It is easy to compute the number of those $k$, $1\le k< (q^2-1)q/2$, $(k,q)=1$, for which $\pm k$ satisfies the conditions (a) and (b) of  Theorem \ref{t2}. For instance, if $q=17$ this number is 192,
and for $q=48$ it is 18424. If we compute the said normalized \DED sums $S(a,b)$ successively, we see that the residue classes $\OV k$ modulo $(q^2-1)q$ appear with different frequencies.
Some of them occur at an early stage of the computation and then repeatedly, others only after many attempts. This has the effect that a long list of pairs $(t,a_1)$ may be necessary for verifying the occurrence of all possible residue classes.
In the case $q=24$, we needed about 65000 normalized \DED sums to verify the occurrence of all 2300 classes, in the case $q=60$, however, about 2.5 million normalized \DED sums for all 28792 classes
(in all cases up to the factor $\pm 1$).

\MN
{\em Example.} Let $q=7$, hence $(q^2-1)q=336$. In this case there are 12 numbers $k, 1\le k<336/2$, $(k,7)=1$, such that $\pm k$ satisfies conditions (a) and (b) of Theorem \ref{t2}, namely
$k=6,18,30,54,66, 78, 90, 102, 114, 138, 150, 162$. In the following table we list the smallest numbers $t$ together with $a$ and $b$, for which the normalized \DED sum $S(a,b)=k'/7$
is such that $k'\equiv \pm k$ mod 336. We also note which of the signs $\pm$ holds here. Recall that $a$ has the form $a=a_1+tj$, $j=0,\LD,6$, $(a,7)=1$, where $a_1^2+1\equiv 0\mod t$. Further, $b=7(a^2+1)/t$.
The table shows that all numbers $k$ appear with $t=1,2, 5$.

\bigskip
\begin{tabular}{r|r|r|r|r}

$t$ &  $a$  & $b$  & $S(a,b)$  &  $k$ \\
\hline
 1 & 1      & 14   & 78/7       & 78   \rule{0mm}{5mm}\\
 1 & 2      & 35   & 102/7     & 102  \rule{0mm}{5mm}\\
 1 & 3      & 70   & 138/7     & 138  \rule{0mm}{5mm}\\
 2  &  1    &  7   & 30/7      &  30 \rule{0mm}{5mm}\\
 2  &  3    & 35   & 66/7      & 66  \rule{0mm}{5mm}\\
 2  &  5    & 91   & 90/7      & 90  \rule{0mm}{5mm}\\
 5  &  2    &  7   &  6/7      & 6   \rule{0mm}{5mm}\\
 5  & 12    & 203  & 162/7     & 162 \rule{0mm}{5mm}\\
 5  & 17    & 406  & 186/7     & $-$150 \rule{0mm}{5mm}\\
 5  & 22    & 679  & 222/7     & $-$114 \rule{0mm}{5mm}\\
 5  & 27    & 1022 & 282/7     & $-$54  \rule{0mm}{5mm}\\
 5  & 32    & 1435 & 318/7     & $-$18  \rule{0mm}{5mm}\\
\end{tabular}


\vspace{0.5cm}
\noindent
Kurt Girstmair            \\
Institut f\"ur Mathematik \\
Universit\"at Innsbruck   \\
Technikerstr. 13/7        \\
A-6020 Innsbruck, Austria \\
Kurt.Girstmair@uibk.ac.at

\end{document}